\documentclass{amsart}
\usepackage{amsmath}
\usepackage{amscd}

\newtheorem{theorem}{Theorem}[section]
\newtheorem{corollary}[theorem]{Corollary}
\newtheorem{lemma}[theorem]{Lemma}
\newtheorem{proposition}[theorem]{Proposition}
\newtheorem{definition}[theorem]{Definition}

\newtheorem{remark}[theorem]{Remark}

\newcommand{\R}{\mathbb{R}}
\newcommand{\Z}{\mathbb{Z}}

\newcommand{\N}{\mathbb{N}}

\newcommand{\Hb}{\mathcal{H}}

\begin{document}


\title[Coarse embeddings of metric spaces into Banach spaces]
{\Large Coarse embeddings of metric spaces into Banach spaces}

\author{Piotr W. Nowak}

\address{Institute of Mathematics, Warsaw University, Banacha 2, 02-097 Warsaw, Poland.}

\address{Department of Mathematics,
Tulane University, 6823 St. Charles Avenue, New Orleans, LA 70118,
USA.}

\email{pnowak@math.tulane.edu}



\begin{abstract}
There are several characterizations of coarse embeddability of a
discrete metric space into a Hilbert space. In this note we give
such characterizations for general metric spaces. By applying
these results to the spaces $L_p(\mu)$, we get their coarse
embeddability into a Hilbert space for $0<p<2$. This together with
a theorem by Banach and Mazur yields that coarse embeddability
into $\ell_2$ and into $L_p(0,1)$ are equivalent when $1 \le p<2$.
A theorem by G.Yu and the above allow to extend to $L_p(\mu)$,
$0<p<2$, the range of spaces, coarse embedding into which
guarantees for a finitely generated group $\Gamma$ to satisfy the
Novikov Conjecture.
\end{abstract}
\maketitle

\section{Introduction}
Coarse embeddings were introduced by M.~Gromov.
\begin{definition}[{\cite[7.E.]{gromov}}]\label{defnition - coarse embedding}
\normalfont Let $X,Y$ be metric spaces. A (not necessarily continuous function)
$f\colon X\to Y$ is a \emph{coarse embedding} if there exist non-decreasing
functions $\rho_1,\rho_2\colon [0,\infty)\to[0,\infty)$ satisfying
\begin{enumerate}
\item $\rho_1(d_X(x,y))\le d_Y(f(x),f(y))\le \rho_2(d_X(x,y))$\; for all \;$x,y\in X$,
\item $\lim_{t\to\infty}\rho_1(t)=+\infty$.
\end{enumerate}

\end{definition}
In the language of coarse geometry, a coarse embedding $f\colon X\to Y$
is a coarse equivalence of $X$ and $f(X)$ (for the notions of coarse geometry
we refer the reader to \cite{gromov,roe}).

Coarse embeddings have been in the spotlight recently due
to a remarkable theorem by G.~Yu \cite{yu}, who showed
that every discrete metric space $\Gamma$, which embeds coarsely into a
Hilbert space, satisfies the Coarse Baum-Connes Conjecture.
In particular, if $\Gamma$ is a finitely generated group with word length metric
coarsely embeddable into a Hilbert space and the classifying space $B\Gamma$ has
a homotopy type of a finite CW-complex, then
the Novikov Conjecture holds for $\Gamma$.

Gromov asked \cite[7.E$_2$.]{gromov} if every separable metric
space embeds coarsely into a Hilbert space. A negative answer to
this question was given by A.N.~Dranishnikov, G.~Gong,
V.~Lafforgue and G.~Yu \cite{dgly}. The counterexample was
constructed using ideas of P.~Enflo, who answered negatively a
question raised by J.~Smirnov, whether every separable metric
space is uniformly homeomorphic to a subset of a Hilbert space
\cite{enflo}.

Inspired by the work of Dranishnikov, Gong, Lafforgue and Yu
\cite{dgly}, we go somewhat further in this direction and try to
find similarities between coarse embeddings and uniform embeddings
(i.e., uniform homeomorphisms onto a subset). As a result, we
generalize characterizations of coarse embeddability into $\ell_2$
given in \cite{dgly} in terms of negative definite functions from locally finite to general metric space.
In particular we show, that coarse embeddability of an arbitrary
metric space depends only on its finite subset. As an application
we prove that for $0<p<2$ the space $L_p$ admits a coarse
embedding into a Hilbert space. This and the fact that the space $L_p(0,1)$ for $p>1$ contains $\ell_2$ as a subspace
yields that coarse embeddability into $\ell_2$ and into $L_p(0,1)$ are equivalent.
All of the above combined with Yu's theorem on the Baum-Connes Conjecture \cite{yu} and the Descent Principle
\cite{roe} implies obviously that the Novikov Conjecture holds
for finitely generated groups coarsely embeddable into a
separable $L_p(\mu)$-space, $0<p<2$.

Throughout this note we will talk about \emph{coarse
embeddability}, having in mind that the embedding is into a
Hilbert space. We will also assume, that all Hilbert spaces
are real - this makes no loss to generality, since embedding into a complex Hilbert space
is equivalent to embedding into a real Hilbert space.

Let us also note that if $f\colon X\to\Hb$ is a coarse embedding
then a reasonable assumption on $X$, allowing $\Hb=\ell_2$, is for
$X$ to contain an at most countable $c$-net $\mathcal{C}$ for some
constant $c>0$ (i.e., for every $x\in X$ there exists
$y\in\mathcal{C}$ such that $d(x,y)\le c\;$). Of course, in the
coarse category $\mathcal{C}$ and $X$ are isomorphic. In
particular, this is the case, when our metric space $X$ is a
(finitely generated) group $\Gamma$ with a word length metric or
(more generally) when $X$ is separable.
\\\\\\
\textbf{Acknowledgements.} This note is a part of the author's
M.Sc. $\!\!\!$ thesis, written
under supervision of prof. Henryk Toru\'{n}czyk
at the Institute of Mathematics of the Warsaw University in Poland and defended in June, 2003.\\
The author wishes to express his gratitude to prof. Toru\'{n}czyk for
guidance and an uncountable number of valuable sugesstions and remarks. The author is also very grateful to
prof. Tadeusz Ko\'{z}niewski and prof. Zbigniew Marciniak for
their active interest in the publication of this paper
and to prof. Guoliang Yu for encouragement.\\

\section{Preliminaries}
Positive definite functions appeared in the work of E.H.~Moore \cite{moore}.
Let us briefly recall basic definitions and
facts concerning positive and negative definite kernels and
functions. For more details we refer the reader to
\cite{schoenberg2}, \cite[Chapter 8]{bl} and \cite{amm}.

\begin{definition}\normalfont
By a kernel on a set $X$ we
mean any symmetric function \linebreak $K:X\times X\to\R$.
A kernel $K$ is said to be:
\begin{enumerate}
\item \emph{positive definite} if
$\sum K(x_i,x_j)c_ic_j\ge 0$ for all $n\in\N$ and $x_1,...,x_n\in X$,
$c_1,...,c_n\in \R$.
\item \emph{negative definite} if
$\sum K(x_i,x_j)c_ic_j\le 0$ for all $n\in\N$ and $x_1,...,x_n\in X$,
$c_1,...,c_n\in \R$ such that $\sum c_i=0$.
\end{enumerate}
A function $f\colon X\to \R$ on a commutative metric group $X$
is said to be positive (negative)
definite if $K(x,y)=f(x-y)$ is a positive
(negative) definite kernel.
\end{definition}

A positive (negative) definite kernel $K$ is said to be
\emph{normalized} if $K(x,x)=1$ ($K(x,x)=0$) for all $x\in X$. The
set of positive definite kernels is closed under limits
in the function space $\R^{X\times X}$ with pointwise convergence.\\

The following theorem describes the relation between negative and
positive definite kernels.

\begin{theorem}[I.J.~Schoenberg \cite{schoenberg2}; see also {\cite[Proposition 8.4]{bl}}]\label{characterization of neg.def.kernel by exp}
A kernel $N$ on $X$ is negative definite  if and only if $e^{-tN}$
is positive definite for every $t>0$.
\end{theorem}

The key to study the connection between maps into Hilbert spaces
and positive and negative definite kernels is the following
theorem (see \cite[Proposition 8.5]{bl}).

\begin{theorem}[(1) E.H.~Moore \cite{moore}; (2) I.J.~Schoenberg \cite{schoenberg2}]\label{moore-schoenberg theorem}
\textnormal{(1)} A kernel $K$ on a set $X$ is positive definite
if and only if  there exist a Hilbert space $\Hb$ and a map
$f\colon X\to\Hb$ such that
$$K(x,y)=\langle f(x),f(y)\rangle\;\;\;\;\text{ for all } x,y\in X.$$
\textnormal{(2)} A normalized kernel $N$ on $X$ is negative
definite if and only if there exist a Hilbert space $\Hb$
and a map $f\colon X\to\Hb$ such that
$$N(x,y)=\Vert f(x)-f(y)\Vert^2\;\;\;\;\text{ for all } x,y\in X.$$
\end{theorem}

Part (1) of Theorem \ref{moore-schoenberg theorem} is the main
tool in the study of uniform embeddings (i.e., uniform
homeomorphisms onto a subset) into Hilbert spaces (see \cite{amm,
bl}) while part (2) gives, as an immediate consequence, the
well-known Schoenberg's characterization of metric spaces which
embed isometrically into a Hilbert space: \emph{a metric space
$(X,d)$ embeds isometrically into a Hilbert space if and only if
$d^2$ is a negative definite kernel on $X$.}
\newline

\section{Coarse embeddings into Hilbert spaces \\and negative definite kernels}

In this section we characterize coarse embeddability of metric
spaces. In \cite{dgly} the authors give a condition equivalent to
coarse embeddability of a locally finite metric space into a
Hilbert space in terms of negative definite kernels \cite[Theorem
2.2]{dgly} (recall that a metric space is locally finite if every
ball has finitely many elements). This characterization remains true for
general metric spaces.

\begin{proposition}\label{embeddings first characterization}
A metric space $X$ admits a coarse embedding into a Hilbert space
if and only if there exist a normalized, negative definite kernel
on $X$ and non-decreasing functions $\rho_i\colon [0,\infty)\to
[0,\infty)$, $i=1,2$, satisfying
\begin{enumerate}
\item $\rho_1(d(x,y))\le N(x,y)\le \rho_2(d(x,y))$; \item
$\lim_{t\to\infty}\rho_1(t)=\infty$.
\end{enumerate}
\end{proposition}
\begin{proof} Suppose, that there exists a kernel $N$ satisfying conditions (1)
and (2). By Theorem \ref{moore-schoenberg theorem}, there exist a
Hilbert space $\Hb$ and a map $f\colon X\to\Hb$ such that
$$N(x,y)=\Vert f(x)-f(y)\Vert^2$$
for all $x,y\in X$. In other words, $f$ satisfies
\newline
$$\sqrt{\rho_1(d(x,y))}\le\Vert f(x)-f(y)\Vert \le \sqrt{\rho_2(d(x,y))}.$$
\newline
Obviously, $f$ is a coarse embedding.

Now suppose, that $f\colon X\to\Hb$ is a coarse embedding and define
$N(x,y)=\Vert  f(x)-f(y)\Vert^2$. Then, $N$ is a normalized,
negative definite kernel on $X$: when $\sum c_i=0$ and $x_1,...,x_n\in X$
then

\begin{eqnarray*}
\sum_{i,j=1}^n \Vert f(x_i)\Vert^2 c_ic_j&=&
\left(\sum_{i=1}^n c_i\right)\left(
\sum_{j=1}^n\Vert f(x_i)\Vert^2 c_j\right),
\end{eqnarray*}
hence
\begin{eqnarray*}
\sum_{i,j=1}^n  N(x_i,x_j)c_i c_j=\\
\sum_{i,j=1}^n & \left( \Vert f(x_i)\Vert^2+\Vert f(x_j)\Vert^2-2
\langle f(x_i),f(x_j)\rangle\right)c_ic_j=\\
&=-2\left\Vert \sum\limits_{i=1}^n c_i f(x_i)\right\Vert^2\le 0
\end{eqnarray*}
\end{proof}

\begin{remark}\normalfont
The proof of Theorem \ref{moore-schoenberg theorem} and Theorem
2.2 \cite{dgly} use the same reasoning. In fact, the proof of the
latter works also in case of general metric spaces, except that
for not locally finite spaces the construction may result in a
non-separable Hilbert space.
\end{remark}

Proposition \ref{embeddings first characterization} leads to the
following
\begin{corollary}\label{cor - metric neg def}
Let $(X,d)$ be a metric space. If for some non-decreasing function
$\alpha:[0,\infty)\to [0,\infty)$ satisfying $\lim_{t\to\infty} \alpha(t)=\infty$
the kernel \ $\alpha\circ d$ is negative definite, then $X$ admits a coarse embedding
into a Hilbert space.
\end{corollary}
\begin{proof}
If $\alpha\circ d$ is a negative kernel on $X$, then it satisfies the
conditions of Proposition \ref{embeddings first characterization} with
$\rho_1(t)=\rho_2(t)= \alpha(t)$.
\end{proof}

Also in \cite{dgly} the authors prove that coarse embeddability
in the case of locally finite metric spaces can be reduced to the
question, whether all finite subsets can be mapped into a Hilbert space
by functions having uniform estimates on the distance of images of two
points \cite[Definition 3.1, Proposition 3.2]{dgly}.

We will show that again the local finiteness assumption can be
dropped and that coarse embeddability of a metric space
into a Hilbert space is determined by its finite subsets.
\begin{theorem}\label{coarse-embeddability-finite-subsets}
A metric space $X$ admits a coarse embedding into a Hilbert space
if and only if there exist non-decreasing functions
$\rho_i\colon [0,\infty)\to [0,\infty)$, $i=1,2$, such that
 $\;\lim_{t\to\infty}\rho_1(t)=\infty$ and for every finite subset
$A\subset X$ there exists a map $f_A\colon A\to\ell_2$ satisfying
$$ \rho_1(d(x,y))\le \Vert f_A(x)-f_A(y)\Vert \le \rho_2(d(x,y))$$
for every $x,y\in X$.
\end{theorem}
\begin{proof}
Suppose that there exist functions $\rho_1, \rho_2$ and
$f_A\colon A\to\ell_2$ as above. The kernel $N_A(x,y)=\Vert
f_A(x)-f_A(y)\Vert^2$ is negative definite on $A$ and we have
$$ \rho^2_1(d(x,y))\le N_A(x,y) \le \rho^2_2(d(x,y))$$
for all $x,y\in A$. Define a new kernel
$K_A=e^{-N_A}$. By Theorem \ref{characterization of neg.def.kernel by exp},
the kernel $K_A$ is positive definite and
$$e^{-\rho^2_2(d(x,y))}\le K_A (x,y) \le e^{-\rho^2_1(d(x,y))}.$$
Extend the kernel $K_A$ to a positive definite kernel
$\widetilde{K}_A$ on $X$, by defining it to be $0$ outside the set
$A\times A$. For all $x,y\in X$ we have $\vert
\widetilde{K}_A(x,y)\vert\le 1$, thus the closure of the set
$\left\{ \widetilde{K}_A\colon  \#A<\infty\right\}$ in the space
$\R^{X\times X}$ with pointwise convergence topology is compact,
by Tichonov's theorem. Let $K$ be a pointwise limit of a
convergent subnet of the net\footnote{Here, by a \emph{net} we mean a generalized sequence}
$\left\{\widetilde{K}_A\colon
\#A<\infty\right\}$. Then $K$ is a positive definite kernel on
$X$, satisfying
$$ e^{-\rho^2_2(d(x,y))}\le K(x,y) \le e^{-\rho^2_1(d(x,y))}$$
for all $x,y\in X$.\\
Since $K(x,y)>0$ for all $x,y\in X$, we can set $N(x,y)=-\ln
K(x,y)$. We will show that $N$ is negative definite. Indeed, since
$K$ is a pointwise limit and the function $e^{-x}$ is continuous,
for every $t>0$ we have

\begin{equation}\label{net-limit-equation}
e^{-tN(x,y)}=K(x,y)^t=\left(\,\lim
\widetilde{K}_A(x,y)\right)^t
=\lim \left(\widetilde{K}_A(x,y)\right)^t\\
\end{equation}
where the limit is the one of the convergent subnet. By Theorem
\ref{characterization of neg.def.kernel by exp}
the kernels $\widetilde{K}_A^{\,t}$ are positive definite for every $t>0$ and the
right-hand side of equation (\ref{net-limit-equation})
is a positive definite kernel as a pointwise limit
of positive definite kernels .
Thus $e^{-tN(x,y)}$ is positive
definite for every $t>0$ and Theorem \ref{characterization of
neg.def.kernel by exp} yields that $N$ is a negative definite
kernel on $X$ satisfying
$$\rho^2_1(d(x,y))\le N(x,y) \le \rho^2_2(d(x,y)).$$
By Theorem \ref{embeddings first characterization}, $X$ admits
a coarse embedding into a Hilbert space.\\
The other implication is obvious.
\end{proof}

\section{Coarse embeddings of $L_p$-spaces}\label{section-banach-spaces}
In this section, as an application of the above results, we show
that the space $L_p(\mu)$ is coarsely embeddable for $0<p\le 2$, which
gives, as an immediate consequence, a stronger version of
\cite[Proposition 3.3]{dgly},  where the authors showed, that
every locally finite subspace of $\ell_1(\N)$ admits a coarse
embedding into a Hilbert space.

The spaces $L_p(\mu)$ are metric only for $p\ge 1$ but we will consider also the
case $0<p<1$.
\begin{proposition}\label{theorem on Lp spaces}
For any measure $\mu$, the space $L_p{(\mu)}$ (in particular
$\ell_p$) admits a coarse embedding into a Hilbert space for\;
$0< p<2$.
\end{proposition}
To prove Proposition \ref{theorem on Lp spaces} we need the following lemma, which
was first proved by I.J.~Schoenberg \cite[Theorem 2]{schoenberg1}
\begin{lemma}
Let $N$ be a negative definite kernel on $X$ and $N(x,y)\ge 0$ for
all $x,y\in X$. Then the kernel $N^{\alpha}$ is negative definite
for any $0<\alpha<1$.
\end{lemma}
\begin{proof}[Proof of the Lemma]
Let $N$ be a negative definite kernel.
 Then for every $t\ge 0$ the kernel $1-e^{-tN}\ge 0$ is also
negative definite, and we have
$$\int\limits_0^{\infty}{\left(1-e^{-tN}\right)d\mu(t)}\ge 0$$
for every positive measure $\mu$ on $[0,\infty)$. For every $x>0$
and $0<\alpha<1$ the following formula holds
$$x^{\alpha}=c_{\alpha}\int\limits_0^{\infty}{\left(1-e^{-tx}\right)t^{-\alpha-1} dt},$$
where $c_{\alpha}$ is some positive constant. Thus $N^{\alpha}$ is
also a negative definite kernel for every $0<\alpha<1$.
\end{proof}

\begin{proof}[Proof of Proposition \ref{theorem on Lp spaces}]
The kernel $\vert x-y\vert^2$ is negative definite on the real
line (as a square of the metric on a Hilbert space). By the above
Lemma, for any $0<p\le 2$ the kernel $\vert x-y\vert^p$ also is
negative definite on $\R$, i.e.,
$$\sum \vert x_i-x_j\vert^p c_i c_j\le 0$$
for every such $p$, all $x_1,...,x_n\in \R$ and $c_1,...,c_n\in\R$
such that $\sum c_i=0$. By integrating the above inequality with
respect to measure $\mu$, we get that the function $\Vert
x\Vert^p$ is negative definite on $L_p(\mu)$ and the assertion
follows from Corollary \ref{cor - metric neg def}.
\end{proof}

\begin{corollary}\label{corollary-on-Lp-spaces}
If a metric space $X$ admits a coarse embedding into the space
$L_p(\mu)$ for $0<p< 2$ and some measure $\mu$ (in particular
into $\ell_p$), then $X$ admits a coarse embedding into a Hilbert
space.
\end{corollary}

G.Yu showed that every discrete metric space $\Gamma$ with bounded geometry
which admits a coarse embedding into a Hilbert space satisfies
the Coarse Baum-Connes Conjecture \cite{yu}. In particular, if $\Gamma$
is a finitely generated group with word length metric then the
Descent Principle \cite{roe} implies that the
Novikov Conjecture holds for $\Gamma$.  To situate the results of
this section in the context of the Novikov Conjecture we give this
obvious
\begin{corollary}
Let $\Gamma$ be a finitely generated group and let the classifying
space $B\Gamma$ have the homotopy type of a finite CW-complex. If $\Gamma$ (as a metric space with word
length metric) admits a coarse embedding into a (separable) space
$L_p(\mu)$, $0<p< 2$
then the Novikov Conjecture holds for $\Gamma$.
\end{corollary}

Let us also note an interesting consequence of a theorem due to
S.Banach and S.Mazur on linear dimension of the spaces $L_p$
\cite[Theorem 6, p.123]{banach} which states that $\ell_2$ is
isomorphic to a subspace $L_p(0,1)$ if $p>1$. From this we
conclude, that every separable metric space, coarsely embeddable
into a Hilbert space, admits a coarse embedding into any
$L_p(0,1)$ with $p\ge 1$. Moreover, for $1<
p\le 2$ the space $L_p(0,1)$ is isometric to a subspace of
$L_1(0,1)$ (\cite[Theorem II.3.14, p.139]{lindenstrauss-tzafriri}).
Thus, we arrive at the following

\begin{corollary}
Let $X$ be a separable metric space. The following statements are
equivalent:
\begin{enumerate}
\item $X$ is coarsely embeddable into $\ell_2$;
\item $X$ is coarsely embeddable into $L_p(0,1)$ for some (equivalently all) $1\le p<2$.
\end{enumerate}

\end{corollary}

\section{Final remarks}
\textbf{A quick example of a discrete metric space not coarsely
embeddable into $\ell_2$.} Let $\mathcal{C}$
denote the set of all integer sequences which are eventually zero.
This is the abelian group $\oplus_{n\in\N}\Z$ on which the word length metric (with
respect to the standard generators) is just the $\ell_1$ metric.
With this metric $\mathcal{C}$ is coarsely embeddable into
$\ell_2$ by Proposition \ref{theorem on Lp spaces}. There is another
natural metric on $\mathcal{C}$, induced by the $\sup$-norm
$\Vert\cdot\Vert_0$. We claim that
$(\mathcal{C},\Vert\cdot\Vert_0)$ is not coarsely embeddable into
$\ell_2$.

To show this consider $(\mathcal{C},\Vert\cdot\Vert_0)$ as a 1-net in the
space $c_0$. If $(\mathcal{C},\Vert\cdot\Vert_0)$ would be coarsely embeddable
into $\ell_2$ then so would be $c_0$. By a result of I.~Aharoni \cite{aharoni}
there exists a constant $K>0$ such that for every separable metric space $X$ there is
a map $T\colon X\to c_0$ satisfying
$$d(x,y)\le \Vert Tx-Ty\Vert\le K d(x,y) \ \ \ \ \text{for all } x,y\in X.$$
Therefore any separable metric space would be coarsely embeddable into $\ell_2$,
which is not the case by a result in \cite{dgly}. This contradiction completes
the proof.

Let us note that, unfortunately, the group $\mathcal{C}$ is infinitely generated
and fails to be locally finite.\newline

(A) The results in Section \ref{section-banach-spaces} are
counterparts of facts from the theory of "classical" uniform
embeddings into Hilbert space (see \cite{amm, bl}). A theorem of
I.~Aharoni, B.~Maurey and B.S.~Mityagin  \cite[Corollary 4.3]{amm}
states, that $L_p(\mu)$ is not uniformly homeomorphic to any
subset of a Hilbert space if $p>2$. It is natural to ask, whether
this result has a counterpart for coarse uniform embeddings, i.e.,
does the space $L_p(\mu)$ admit a coarse uniform embedding into a
Hilbert space if $p>2$?\newline

(B) Theorem \ref{embeddings first characterization} can
provide an alternative proof of the characterization of coarse
embeddability, given by M.~Dadarlat and E.~Guentner
\cite[Proposition 2.1]{dg}: \emph{a metric space $X$ is coarsely
embeddable into a Hilbert space $\Hb$ if and only if for every
$R>0$ and every $\varepsilon>0$ there exist $S>0$ and a map
$f\colon X\to\Hb$ such, that for all $x,y\in X$ we have $\Vert
f(x)\Vert=1$, $\sup\{\,\Vert f(x)-f(y)\Vert\colon d(x,y)\le
R\,\}\le\varepsilon$ and
$\lim_{S\to\infty} \inf \{\,\Vert f(x)-f(y)\Vert\colon d(x,y)\ge S\}=\sqrt{2}$.}
We leave the details to the reader.\newline

(C) Since $c_0^{\,*}=\ell_1$ and $\ell_1^{\,*}=\ell_{\infty}\supset c_0$,
Proposition \ref{theorem on Lp spaces} and the above discussion
 provide an example of a Banach space, which is
coarsely embeddable, while it's dual is not, and of a Banach
space, which is not coarsely embeddable, while it's dual
is.\newline

\bibliographystyle{amsalpha}

\end{document}